\documentclass[a4paper,conference]{IEEEtran}
\usepackage{amssymb}
\usepackage[cmex10]{amsmath}
\newtheorem{theorem}{Theorem}
\begin{document}
\title{Kernels of Integral Equations Can Be \\
Boundedly Infinitely Differentiable on $\mathbb{R}^2$}
\author{\IEEEauthorblockN{Igor M. Novitskii}
\IEEEauthorblockA
{Khabarovsk Division \\
Institute of Applied Mathematics \\
Far-Eastern Branch of the Russian Academy of Sciences \\
Khabarovsk, Russia\\
Email: novim@iam.khv.ru}
}
\maketitle
\begin{abstract}
\boldmath
In this paper, we reduce the general linear integral equation of the third kind 
in $L^2(Y,\mu)$, with largely arbitrary kernel and coefficient, to an equivalent 
integral equation either of the second kind or of the first kind in 
$L^2(\mathbb{R})$, with the kernel being the linear pencil of bounded infinitely
differentiable bi-Carleman kernels expandable in absolutely and uniformly 
convergent bilinear series. The reduction is done by using unitary equivalence 
transformations.
\end{abstract}
\IEEEpeerreviewmaketitle
\section{Introduction}
In the theory of general linear integral equations in $L^2$ spaces, equations 
with bounded infinitely differentiable bi-Carleman kernels (termed $K^\infty$ 
kernels) should and do lend themselves well to solution by approximation and 
variational methods. The question of whether a second-kind integral equation 
with arbitrary kernel can be reduced to an equivalent one with a $K^\infty$ 
kernel was positively answered using a unitary-reduction method by the author 
\cite{nov:IJPAM2}. In the present paper, our goal is to extend the method in
order to deal with a third-kind integral equation (\eqref{Eq1} below) with
arbitrary measurable kernel and coefficient.
\par
Results obtained are presented with proofs in Section~\ref{RedTHS}.
The results say that the general linear integral equation of the third kind
in $L^2(Y,\mu)$  can be reduced to an equivalent integral equation  either of 
the second kind (Theorem~\ref{thirdkind}) or of the first kind (Theorem~\ref{firstkind})
in $L^2(\mathbb{R})$, with the kernel being the linear pencil of $K^\infty$ kernels
of Mercer type or of Hilbert-Schmidt $K^\infty$ kernels of Mercer type,
respectively.
\par
Before we can write down and prove our results, we need to fix the terminology 
and notation and to give some definitions and preliminary material.
\section{Preliminary Notions and Results}
\subsection{Spaces}
Throughout this paper, $\mathcal{H}$ is a complex, separable, infinite-dimensional
Hilbert space with norm $\|\cdot\|_{\mathcal{H}}$ and inner product
$\left\langle\cdot,\cdot\right\rangle_{\mathcal{H}}$.
$(Y,\mu)$ is a measure space $Y$ equipped with a positive,
$\sigma$-finite, complete, separable, and nonatomic, measure $\mu$.
$L^2(Y,\mu)$ is the Hilbert space of (equivalence classes of)
$\mu$-measurable complex-valued functions on $Y$ equipped with the inner product
$
\left\langle f,g\right\rangle_{L^2(Y,\mu)}=\int_Y f(y)\overline{g(y)}\,d\mu(y)
$
and the norm
$\left\|f\right\|_{L^2(Y,\mu)}=\left\langle f,f\right\rangle_{L^2(Y,\mu)}^{1/2}$;
when $\mu$ is the Lebesgue measure on the real line $\mathbb{R}$,
$L^2(\mathbb{R},\mu)$ is abbreviated into $L^2$, and $d\mu(y)$ into $dy$.
$C(X,B)$, where $B$ is a Banach space with norm
$\|\cdot\|_B$, is the Banach space (with the norm
$\|f\|_{C(X,B)}=\sup_{x\in X}\,\|f(x)\|_B$) of continuous $B$-valued functions
defined on a locally compact space $X$ and \textit{vanishing at infinity}
(that is, given any $f\in C(X,B)$ and $\varepsilon>0$, there exists a compact
subset $X(\varepsilon,f) \subset X$ such that $\|f(x)\|_{B}<\varepsilon$
whenever $x\not\in X(\varepsilon,f)$). A series $\sum_n f_n$ is \textit{$B$-absolutely convergent in $C(X,B)$}
if $f_n\in C(X,B)$ ($n\in\mathbb{N}$) and the series $\sum_n \|f_n(x)\|_B$
converges in $C(X,\mathbb{R})$. Given an equivalence class $f\in L^2$ containing a function of
$C(\mathbb{R},\mathbb{C})$, the symbols $[f]$ and $[f]^{(i)}$ are used to denote
that function and its $i$th derivative, if exists. The symbols $\mathbb{C}$
and $\mathbb{N}$ refer to the complex plane and the set of all positive
integers, respectively.
\subsection{Linear Operators}
Throughout, $\mathfrak{R}(\mathcal{H})$ denotes the Banach algebra of bounded
linear operators acting on $\mathcal{H}$. For an operator $T$ of
$\mathfrak{R}(\mathcal{H})$, $T^*$ stands for the adjoint to $T$ (w.r.t.
$\langle\cdot,\cdot\rangle_{\mathcal{H}}$), and the family
$\mathcal{M}^{+}(T)$ is defined as the set of all those operators
$P\in\mathfrak{R}(\mathcal{H})$ that are positive (that is,
$\langle Px,x\rangle_{\mathcal{H}}\ge 0$ for all $x\in\mathcal{H}$) and
representable in one or other of the forms $P=TB$, $P=BT$, where
$B\in\mathfrak{R}(\mathcal{H})$. A factorization of an operator
$T\in\mathfrak{R}({\mathcal{H}})$ into the product $T = WV^*$ ($V$,
$W\in \mathfrak{R}({\mathcal{H}})$) is called an \textit{$\mathcal{M}$
factorization} for $T$  provided that $VV^*$,
$WW^*\in\mathcal{M}^{+}(T)$ \cite{nov:IeJPAM}. An ingenuous example of an
$\mathcal{M}$ factorization for any $T\in\mathfrak{R}(\mathcal{H})$ is given
by taking $W=UP$, $V=P$, where $P$ is the positive square root of
$|T|=(T^*T)^{\frac1{2}}$ and $U$ is the partially isometric factor in the
polar decomposition $T=U|T|$.
\par
A bounded linear operator $U\colon\mathcal{H}\to L^2$ is \textit{unitary} 
if it has range $L^2$ and 
$\left\langle Uf,Ug\right\rangle_{L^2}=\left\langle f,g\right\rangle_{\mathcal{H}}$
for all $f$, $g\in\mathcal{H}$. An operator $S\in \mathfrak{R}({\mathcal{H}})$ is 
\textit{unitarily equivalent} to an operator 
$T\in \mathfrak{R}\left(L^2\right)$ if a unitary operator
$U:{\mathcal{H}}\to L^2$ exists such that $T=USU^{-1}$.
\par
A linear operator $T\colon L^2(Y,\mu)\to L^2(Y,\mu)$ is \textit{integral}
if there is a complex-valued $\mu\times\mu$-measurable
function $\boldsymbol{T}$ (\textit{kernel}) defined on the Cartesian product
$Y^2=Y\times Y$ such that
\begin{equation*}
(Tf)(x)=\int_Y \boldsymbol{T}(x,y)f(y)\,d\mu(y)
\end{equation*}
for every $f$ in $L^2(Y,\mu)$ and $\mu$-almost every $x$ in $Y$ \cite{Halmos:Sun}, \cite{Kor:book1}.
The integral operator $T$ is \textit{bi-integral} if its adjoint $T^*$ is also
an integral operator on $L^2(Y,\mu)$ \cite{Kor:book1}. The integral and the bi-integral operators
are bounded, and need not be compact.
\subsection{Linear Integral Equations}
The linear integral equation of the \textit{third kind} in $L^2(Y,\mu)$
is an equation of the form
\begin{equation}\label{Eq1}
\boldsymbol{H}(x)\phi(x)-\lambda\int_{Y}\boldsymbol{K}(x,y)\phi(y)\,d\mu(y)=\psi(x)\quad
\end{equation}
$\mu$-almost everywhere on $Y$, where  $\boldsymbol{H}\colon Y\to\mathbb{C}$ 
(the \textit{coefficient} of the equation) is a given bounded $\mu$-measurable function,
$\boldsymbol{K}\colon Y^2\to\mathbb{C}$ (the \textit{kernel} of the equation) is a
given kernel of a bi-integral operator $K\in\mathfrak{R}\left(L^2(Y,\mu)\right)$,
the scalar $\lambda\in\mathbb{C}$ (a \textit{parameter}) is given, the function $\psi$ of $L^2(Y,\mu)$
is given, and the function $\phi$ of $L^2(Y,\mu)$ is to be determined.
When the coefficient $\boldsymbol{H}$ has the constant value $0$ (resp., $1$)
$\mu$-almost everywhere on $Y$, the linear integral equation \eqref{Eq1}
is referred to as of the \textit{first} (resp., \textit{second}) \textit{kind}.
\subsection{$K^\infty$ Kernels of Mercer Type}
A \textit{bi-Carleman kernel} $\boldsymbol{T}$ on $Y^2$ is
one for which
\begin{equation*}
\int_Y|\boldsymbol{T}(x,y)|^2d\mu(y)<\infty,\quad
\int_Y|\boldsymbol{T}(y,x)|^2d\mu(y)<\infty
\end{equation*}
for $\mu$-almost every $x$ in $Y$. A \textit{Hilbert-Schmidt kernel}
$\boldsymbol{\varGamma}$ on $Y^2$ is one for
which
\begin{equation*}
\int_Y\int_Y|\boldsymbol{\varGamma}(x,y)|^2d\mu(y)\,d\mu(x)<\infty.
\end{equation*}
A \textit{$K^\infty$ kernel} $\boldsymbol{T}$ is a bi-Carleman kernel on
$\mathbb{R}^2$, which is subject to the following infinite differentiability
requirements:
\begin{enumerate}
\renewcommand{\labelenumi}{(\roman{enumi})}
\item the function $\boldsymbol{T}$ and all its partial derivatives of all
orders are in $C\left(\mathbb{R}^2,\mathbb{C}\right)$,
\item the \textit{Carleman function} $\boldsymbol{t}\colon\mathbb{R}\to L^2$, defined via
$\boldsymbol{T}$ by $\boldsymbol{t}(s)=\overline{\boldsymbol{T}(s,\cdot)}$,
and its (strong) derivatives $\boldsymbol{t}^{(i)}$ of all orders are in
$C\left(\mathbb{R},L^2\right)$,
\item the \textit{Carleman function}
$\boldsymbol{t}^{\boldsymbol{\prime}}\colon\mathbb{R}\to L^2$,
defined via $\boldsymbol{T}$ by $\boldsymbol{t}^{\boldsymbol{\prime}}(s)
=\boldsymbol{T}(\cdot,s)$, and its (strong) derivatives
$(\boldsymbol{t}^{\boldsymbol{\prime}})^{(j)}$ of all orders are in
$C\left(\mathbb{R},L^2\right)$ \cite{nov:CEJM}, \cite{nov:IeJPAM}.
\end{enumerate}
A $K^\infty$ kernel $\boldsymbol{T}$ is called of \textit{Mercer type}
if it induces an integral operator $T\in\mathfrak{R}(L^2)$, with the property
that any operator belonging to $\mathcal{M}^{+}(T)$ is also an integral
operator induced by a $K^\infty$ kernel. Any $K^\infty$ kernel $\boldsymbol{T}$ of Mercer type, along with all its partial and strong
derivatives, is entirely recoverable from the knowledge of at least one
$\mathcal{M}$ factorization for its associated integral operator $T$, by means
of bilinear series formulae universally applicable on arbitrary orthonormal
bases of $L^2$:
\begin{theorem}\label{mfactor}
Let $T\in\mathfrak{R}\left(L^2\right)$ be an integral operator, with a kernel
$\boldsymbol{T}$ that is a $K^\infty$ kernel of Mercer type.
Then, for any $\mathcal{M}$ factorization $T = WV^*$ for $T$ and for any orthonormal
basis $\{u_n\}$ for $L^2$, the following formulae hold
\allowdisplaybreaks
\begin{gather}
\dfrac{\partial^{i+j}\boldsymbol T}{\partial s^i\partial t^j}(s,t)=
\sum_n\left[Wu_n\right]^{(i)}(s)\overline{\left[Vu_n\right]^{(j)}(t)}, \label{meijT1}
\\
\begin{split}\label{meijkt2}
\boldsymbol{t}^{(i)}(s)=\sum_n\overline{\left[Wu_n\right]^{(i)}(s)}Vu_n,
\\
\left(\boldsymbol{t}^{\boldsymbol{\prime}}\right)^{(j)}(t)=
\sum_n \overline{\left[Vu_n\right]^{(j)}(t)}Wu_n
\end{split}
\end{gather}
for all non-negative integers $i$, $j$ and all $s$, $t\in\mathbb{R}$, where
the series of \eqref{meijT1} converges $\mathbb{C}$-absolutely in
$C\left(\mathbb{R}^2,\mathbb{C}\right)$, and the series of \eqref{meijkt2}
converge in $C\left(\mathbb{R},L^2\right)$.
\end{theorem}
\par
The theorem is proven in \cite{nov:IeJPAM}. It can also be seen as
a generalization of both Mercer's \cite{Mercer} theorem (about absoluteness and uniformity of
convergence of bilinear eigenfunction expansions for continuous
compactly supported kernels of positive, integral operators) and Kadota's
\cite{Kadota1} theorem (about term-by-term differentiability of those
expansions while retaining the absolute and the uniform convergence) to various
other settings (for details see \cite{nov:Lon}, \cite{nov:IeJPAM}).
\subsection{An Integral Representation Theorem}
The main device for the proof of our reduction theorems in the next section is
provided by the following result, which characterizes families incorporating
those operators in $\mathfrak{R}({\mathcal{H}})$ that can be simultaneously
transformed by the same unitary equivalence transformation into bi-integral
operators having as kernels $K^\infty$ kernels of Mercer type:
\begin{theorem}\label{infsmooth2}
Suppose that for an operator family
$\left\{S_\gamma\right\}_{\gamma\in\mathcal{G}}$$\subset\mathfrak{R}({\mathcal{H}})$
with an index set of arbitrary cardinality there exists an orthonormal
sequence $\left\{e_n\right\}$ in $\mathcal{H}$ such that
\begin{equation}\label{1.2}
\lim\limits_{n\to\infty}\sup\limits_{\gamma\in\mathcal{G}}\left\|S_\gamma e_n\right\|_{\mathcal{H}}=0,
\quad
\lim_{n\to\infty}\sup_{\gamma\in\mathcal{G}}\left\|S_\gamma^* e_n\right\|_{\mathcal{H}}=0.
\end{equation}
Then there exists a unitary operator $U:\mathcal{H}\to L^2$ such that all the
operators $T_\gamma=U S_\gamma U^{-1}$ $(\gamma\in\mathcal{G})$ and their
linear combinations are bi-integral operators on $L^2$, whose kernels are
$K^\infty$ kernels of Mercer type.
\end{theorem}
\par
The proof is given in \cite{nov:IeJPAM}. It provides an explicit procedure to
find that unitary operator $U:\mathcal{H}\to L^2$ whose existence the
Theorem~\ref{infsmooth2} asserts. The procedure uses no spectral properties of
the operators $S_\gamma$, other than their joint property imposed in
\eqref{1.2}, to determine the action of $U$ by specifying two orthonormal
bases, of $\mathcal{H}$ and of $L^2$, one of which is meant to be the image by
$U$ of the other; the basis for $L^2$ may be chosen to be the Lemari\'e-Meyer
wavelet basis \cite{Ausch}, \cite{Her:Wei}.
\section{Reduction Theorems}\label{RedTHS}
\begin{theorem}\label{thirdkind}
Suppose that the essential range of the coefficient $\boldsymbol{H}$ in
\eqref{Eq1} contains the point $\alpha\in\mathbb C$, that is,
\begin{equation}\label{essval}
\mu\{y\in Y\colon |\boldsymbol{H}(y)-\alpha|<\varepsilon\}>0
\quad\text{for all $\varepsilon>0$.}
\end{equation}
Then equation \eqref{Eq1} is equivalent (via a unitary operator
$U$ from $L^2(Y,\mu)$ onto $L^2$) to a second-kind integral equation in $L^2$,
of the form
\begin{equation}\label{Eq2}
\alpha f(s)+\int_{\mathbb R}\left(\boldsymbol{T}_0(s,t)-\lambda\boldsymbol{T}(s,t)\right)f(t)\,dt=g(s)
\end{equation}
almost everywhere on $\mathbb{R}$, where the function $f$($=U\phi$) of $L^2$
is to be determined, the function $g$($=U\psi$) of $L^2$ is given, both the
functions $\boldsymbol{T}_0$ and $\boldsymbol{T}$ are $K^\infty$ kernels of
Mercer type, not depending on $\lambda$, and the function
$\boldsymbol{T}_0-\lambda\boldsymbol{T}$ is also a $K^\infty$ kernel of Mercer
type.
\end{theorem}
\begin{IEEEproof}
The proof relies primarily on the following observation by Korotkov
\cite[Corollary~1]{Kor:1987}: If $H$ is the multiplication operator induced on $L^2(Y,\mu)$
by the coefficient $\boldsymbol{H}$, and $I$ is identity operator on
$L^2(Y,\mu)$, then the two-element family
$\{S_1=H-\alpha I, S_2=K\}$ of bounded operators on
$\mathcal{H}=L^2(Y,\mu)$ satisfies the assumptions of Theorem~\ref{infsmooth2}.
The construction of Korotkov's sequence $\{e_n\}$ fulfilling
\eqref{1.2} for this family is likely to be of practical use and
deserves to be expounded in some detail.
\par
If $E\subset Y$ is a $\mu$-measurable set of positive finite measure,
the orthonormal sequence of \textit{generalized} Rademacher functions
with supports in $E$ will be denoted by $\{R_{n,E}\}_{n=1}^\infty$  and
is constructed iteratively through successive bisections of $E$ as follows:
$$
R_{1,E}=\frac1{\sqrt{\mu E}}\left(\chi_{E_1}-\chi_{E_2}\right)
$$
provided $E_1\sqcup E_2=E$ with $\mu E_1=\mu E_2=\frac12\mu E$;
$$
R_{2,E}=\frac1{\sqrt{\mu E}}\left(\chi_{E_{1,1}}-\chi_{E_{1,2}}+
\chi_{E_{2,1}}-\chi_{E_{2,2}}\right)
$$
provided $E_{i,1}\sqcup E_{i,2}=E_i$ with $\mu E_{i,k}=\frac14\mu E$
for $i,k=1,2$;
\begin{equation*}
\begin{split}
R_{3,E}=\frac1{\sqrt{\mu E}}(&\chi_{E_{1,1,1}}-\chi_{E_{1,1,2}}+\chi_{E_{1,2,1}}-\chi_{E_{1,2,2}}
\\&
+\chi_{E_{2,1,1}}-\chi_{E_{2,1,2}}+\chi_{E_{2,2,1}}-\chi_{E_{2,2,2}})
\end{split}
\end{equation*}
provided $E_{i,k,1}\sqcup E_{i,k,2}=E_{i,k}$ with
$\mu E_{i,k,j}=\frac18\mu E$ for $i,k,j=1,2$;
and so on indefinitely (here $\chi_Z$ denotes the characteristic function
of a set $Z$ and the unions are disjoint). A relevant result due to Korotkov
states that
\begin{equation}\label{krresult}
\lim\limits_{n\to\infty}\|K^*R_{n,E}\|_{L^2(Y,\mu)}=0
\end{equation}
for any integral operator $K\in\mathfrak{R}\left(L^2(Y,\mu)\right)$ and
any $\mu$-measurable $E\subset Y$ with $0<\mu E<\infty$
(see, e.g., the proof of Theorem~3 in \cite{Kor:1986}).
\par
Let $\{Y_n\}_{n=1}^\infty$ be an ascending sequence of sets of positive finite measure,
such that $Y_n\uparrow Y$, let $\{\varepsilon_n\}_{n=1}^\infty$ be a sequence
of positive reals strictly decreasing to zero, and define
$E_n=Y_n\cap\left\{y\in Y\colon \varepsilon_{n+1}<|\boldsymbol{H}(y)-\alpha|\le \varepsilon_n\right\}$
whenever $n\in\mathbb{N}$. Due to the assumption \eqref{essval}, one can always make
the sets $E_n$ to have finite nonzero measures by an appropriate choice of
$Y_n$ and $\varepsilon_n$ ($n\in\mathbb{N}$). Having done so, let
$e_n=R_{k_n,E_n}$, where, for each $n\in\mathbb{N}$, $k_n$ is an index
satisfying
\begin{equation}\label{Ineq1}
\begin{split}
\|&S_2e_n\|_{L^2(Y,\mu)}+\|S_2^*e_n\|_{L^2(Y,\mu)}
\\&=\|KR_{k_n,E_n}\|_{L^2(Y,\mu)}+\|K^*R_{k_n,E_n}\|_{L^2(Y,\mu)}\le1/n
\end{split}
\end{equation}
(cf. \eqref{krresult}).
Since the $E_n$'s are pairwise disjoint, the $e_n$'s form an orthonormal
sequence in $L^2(Y,\mu)$. Moreover, by construction of sets $E_n$,
\begin{equation}\label{Ineq2}
\begin{split}
\|S_1e_n\|&_{L^2(Y,\mu)}^2=\|S_1^*e_n\|_{L^2(Y,\mu)}^2
\\&=\left(1/{\mu E_n}\right)
\int_{E_n}\left|\boldsymbol{H}(y)-\alpha\right|^2\,d\mu(y)
\le \varepsilon_n^2.
\end{split}
\end{equation}
\par
One can now assert from \eqref{Ineq1}, \eqref{Ineq2} that
$\|S_re_n\|_{L^2(Y,\mu)}\to 0$, $\|S_r^*e_n\|_{L^2(Y,\mu)}\to 0$ as
$n\to\infty$, for $r=1,2$. By Theorem~\ref{infsmooth2}, there is then a
unitary operator $U\colon L^2(Y,\mu)\to L^2$ such that the operators
$T_0=US_1U^{-1}$, $T=US_2U^{-1}$ and their linear combinations are bi-integral
operators with $K^\infty$ kernels of Mercer type. This unitary operator can
also be used to transform the integral equation \eqref{Eq1} into an equivalent
integral equation of the form \eqref{Eq2} in such a way that $\boldsymbol{T}_0$,
$\boldsymbol{T}$, and $\boldsymbol{T}_0-\lambda\boldsymbol{T}$, are just those
$K^\infty$ kernels of Mercer type that induce $T_0$, $T$, and $T_0-\lambda T$,
respectively. In operator notation, such a passage from \eqref{Eq1} to
\eqref{Eq2} looks as follows:
$U\psi=U(H-\lambda K)U^{-1}U\phi=U(\alpha I+S_1-\lambda S_2)U^{-1}U\phi
=\alpha f+(T_0-\lambda T)f=g$ where $f=U\phi$, $g=U\psi$.
The theorem is proved.
\end{IEEEproof}
\begin{theorem}\label{firstkind}
If, with the notation and hypotheses of Theorem~\ref{thirdkind}, $\alpha =0$,
then equation \eqref{Eq1} is equivalent to a first-kind integral equation
in $L^2$, of the form
\begin{equation}\label{Eq4}
\int_{\mathbb{R}}\left(\boldsymbol{\varGamma}_0(s,t)-\lambda\boldsymbol{\varGamma}(s,t)\right)f(t)
\,dt=w(s)
\end{equation}
almost everywhere on $\mathbb{R}$, where the function $f$ of $L^2$ is to be
determined, both the functions $\boldsymbol{\varGamma}_0$ and $\boldsymbol{\varGamma}$ are
Hilbert-Schmidt $K^\infty$ kernels of Mercer type, not depending on $\lambda$,
and the function $\boldsymbol{\varGamma}_0-\lambda\boldsymbol{\varGamma}$ is
also a Hilbert-Schmidt $K^\infty$ kernel of Mercer type.
\end{theorem}
\begin{IEEEproof} In this case the equation \eqref{Eq2}, equivalent to \eqref{Eq1}, becomes
\begin{equation}\label{Eq3}
\int_{\mathbb{R}}\left(\boldsymbol{T}_0(s,t)
-\lambda\boldsymbol{T}(s,t)\right)
f(t)\,dt=g(s)
\end{equation}
for almost all $s$ in $\mathbb{R}$. Let $m\in L^2$ be such that $[m]$ is an
infinitely differentiable, positive function all whose derivatives $[m]^{(i)}$
belong to $C(\mathbb{R},\mathbb{R})$, and let $M$ be the multiplication
operator induced on $L^2$ by $m$. Multiply both sides of equation \eqref{Eq3}
by $m$, to recast it into an equivalent equation of the form \eqref{Eq4},
with the same sought-for function $f\in L^2$, the new right side $w=Mg\in L^2$,
and the new kernel $\boldsymbol{\varGamma}_0-\lambda\boldsymbol{\varGamma}$,
where $\boldsymbol{\varGamma}_0(s,t)=[m](s)\boldsymbol{T}_0(s,t)$,
$\boldsymbol{\varGamma}(s,t)=[m](s)\boldsymbol{T}(s,t)$.
It is to be proved that $\boldsymbol{\varGamma}_0$, $\boldsymbol{\varGamma}$
are Hilbert-Schmidt $K^\infty$ kernels. The proof is further given only for
$\boldsymbol{\varGamma}$, as the proof for the other kernel
$\boldsymbol{\varGamma}_0$ is entirely similar. If $\boldsymbol{t}$ is the
associated Carleman function of the $K^\infty$ kernel $\boldsymbol{T}$
(see (ii)), then
  \begin{multline*}
  \int_{\mathbb{R}}\int_{\mathbb{R}}|\boldsymbol{\varGamma}(s,t)|^2\,dt\,ds=
  \int_{\mathbb{R}}m^2(s)\|\boldsymbol{t}(s)\|_{L^2}^2\,ds
  \\\le
  \|\boldsymbol{t}\|_{C\left(\mathbb{R},L^2\right)}^2
  \|m\|_{L^2}^2<\infty,
  \end{multline*}
implying that $\boldsymbol{\varGamma}$ is a Hilbert-Schmidt kernel and hence
induces two Carleman functions $\boldsymbol{\gamma}$,
$\boldsymbol{\gamma}^{\boldsymbol{\prime}}\colon \mathbb{R}\to L^2$ by
$\boldsymbol{\gamma}(s)=\overline{\boldsymbol{\varGamma}(s,\cdot)}$,
$\boldsymbol{\gamma}^{\boldsymbol{\prime}}(t)=
\boldsymbol{\varGamma}(\cdot,t)$. The series representation of $\boldsymbol{T}$ 
(see \eqref{meijT1} for $i=j=0$) gives rise to a series representation of
$\boldsymbol{\varGamma}$, namely, with the notation of Theorem~\ref{mfactor},
\begin{multline*}
\boldsymbol{\varGamma}(s,t)=
\sum_n\left[MWu_n\right](s)\overline{\left[Vu_n\right](t)}
\\=[m](s)\sum_n\left[Wu_n\right](s)\overline{\left[Vu_n\right](t)}
\end{multline*}
for all $(s,t)\in\mathbb{R}^2$. Moreover, for all non-negative integers $i$,
$j$ and all $s$, $t\in\mathbb{R}$, the following formulae hold
\begin{multline*}
\dfrac{\partial^{i+j}\boldsymbol{\varGamma}}{\partial s^i\partial t^j}(s,t)
=\sum_n\left[MWu_n\right]^{(i)}(s)\overline{\left[Vu_n\right]^{(j)}(t)}
\\=\sum_{r=0}^i\binom{i}{r}[m]^{(i-r)}(s)
\left(\sum_n\left[Wu_n\right]^{(r)}(s)\overline{\left[Vu_n\right]^{(j)}(t)}\right),
\end{multline*}
\begin{multline*}
\boldsymbol{\gamma}^{(i)}(s)=\sum_n\overline{\left[MWu_n\right]^{(i)}(s)}Vu_n
\\=\sum_{r=0}^i\binom{i}{r}[m]^{(i-r)}(s)
\left(\sum_n\overline{\left[Wu_n\right]^{(r)}(s)}Vu_n\right),
\end{multline*}
\begin{multline*}
\left(\boldsymbol{\gamma}^{\boldsymbol{\prime}}\right)^{(j)}(t)=
\sum_n \overline{\left[Vu_n\right]^{(j)}(t)}MWu_n
\\=M\left(\sum_n \overline{\left[Vu_n\right]^{(j)}(t)}Wu_n\right),
\end{multline*}
in as much as the bracketed series above converge absolutely in
$C(\mathbb{R}^2,\mathbb{C})$, as regards the first formula, and in
$C(\mathbb{R},L^2)$, as regards the last two formulae, by Theorem~\ref{mfactor}.
Therefore
$$
\dfrac{\partial^{i+j}\boldsymbol{\varGamma}}{\partial s^i\partial t^j} \in
C\left(\mathbb{R}^2,\mathbb{C}\right),\quad \boldsymbol{\gamma}^{(i)},
\left(\boldsymbol{\gamma}^{\boldsymbol{\prime}}\right)^{(j)}\in
C\left(\mathbb{R},L^2\right)
$$
for all non-negative integers $i$, $j$, implying that $\boldsymbol{\varGamma}$ is a $K^\infty$ kernel.
\par
If not all of the $K^\infty$ kernels $\boldsymbol{\varGamma}_0$,
$\boldsymbol{\varGamma}$,
$\boldsymbol{\varGamma}_0-\lambda\boldsymbol{\varGamma}$ are now  of Mercer type,
apply Theorem~\ref{infsmooth2} to the two-element family $\{S_1=MT_0, S_2=MT\}$
of compact operators on $\mathcal{H}=L^2$, the condition \eqref{1.2} for this
family is satisified by any orthonormal sequence $\left\{e_n\right\}$ in $L^2$.
The proof of the theorem is now complete.
\end{IEEEproof}
\section{Conclusion}
In virtue of Theorems~\ref{thirdkind} and \ref{firstkind}, one can confine
one's attention (with no essential loss of generality) to integral equations,
whose kernels are $K^\infty$ kernels of Mercer type, depending linearly on a
parameter. One of the main technical advantages of dealing with such kernels
is that their restrictions to compact rectangles in $\mathbb{R}^2$ are fully
amenable to the methods of the classical theory of ordinary integral equations,
and approximate their original kernels with respect to $C\left(\mathbb{R}^2,\mathbb{C}\right)$ and
$C\left(\mathbb{R},L^2\right)$ norms. This, for instance, can be used directly
\newpage\noindent
to establish an explicit theory of spectral functions for any Hermitian
$K^\infty$ kernel ($\boldsymbol{T}(s,t)=\overline{\boldsymbol{T}(t,s)}$)
by a development essentially the same as the one given by Carleman
\cite{Carl:book} (see also \cite{Trji}, \cite{Costley}, and \cite{Will}, for
further developments). We believe that with regard to $K^\infty$ kernels of
Mercer type this Carleman's line of development can be extended far beyond
the restrictive assumption of a Hermitian kernel.
\par
For the theory of Fredholm determinant and minors, there are some applications in \cite{nov:FEMS},
\cite{nov:FEMJ1}, \cite{nov:FEMJ2}, and \cite{nov:IJPAM1}.

\end{document}